\documentclass{amsart}
\usepackage{amsmath}
\usepackage{amssymb,amsthm}
\usepackage[a4paper]{geometry}
\usepackage{graphicx}
\usepackage{subfig}
\usepackage{fontenc}
\usepackage{multicol}
\usepackage{hyperref}
\usepackage{color}

\makeatletter 

\newcommand{\TeoremaAmbFinalMarcat}[1]{%
  \expandafter\gdef\csname end#1\endcsname{\@endtheorem}}
  
               {\hfill\rule{2.5mm}{2.5mm} \vspace{\parskip} } 
\newtheorem{theorem}{Theorem}[section]

\theoremstyle{definition}
 \TeoremaAmbFinalMarcat{defi}
 \TeoremaAmbFinalMarcat{ex}
\newtheorem{remark}[theorem]{Remark} \TeoremaAmbFinalMarcat{rem}
 \TeoremaAmbFinalMarcat{conv}

\newenvironment{proclama}[1]{
                \par\vspace{\topsep}\noindent{\bf #1}
                \begin{em}}
                {\end{em}\par\vspace{\topsep}}

\newcommand{\A}{\overline{a}}
\newcommand{\B}{\overline{b}}

\newcommand{\wl}{\mathrm{wl}}
\newcommand{\gl}{\mathrm{G}_{A,B,C}}
\newcommand{\FF}[1]{\mathcal{F}_{#1}}
\newcommand{\abc}{_{A,B,C}}
\begin{document}
\title{Experiments suggesting that the distribution of the hyperbolic length of closed geodesics sampling by  word length is Gaussian}

\author{Moira Chas$^1$, Keren Li$^2$ and Bernard Maskit$^3$}
\address{Stony Brook University, Stony Brook, NY 11794\\$^1$ moira@math.sunysb.edu\\$^2$ keren.li@stonybrook.edu\\$^3$Bernard.Maskit@stonybrook.edu}
\date{\today}
\thanks{Supported by NSF grant DMS  - 1105772 }

\maketitle

\begin{abstract}  
Each free homotopy  class of directed closed curves on a surface with boundary can be described by a cyclic reduced word in the generators of the
fundamental group and their inverses. The word  length is the number of letters of the cyclic word.

If the surface  has a hyperbolic metric with geodesic boundary, the  geometric length of the class is the 
length of the unique geodesic.

By computer experiments, we investigate the distribution of  the geometric length  among all classes
with a given word length in the pair of pants surface. Our experiments   strongly suggest that the distribution is normal.\end{abstract}

\section{Introduction and statement of the conjecture}

The fundamental group of a pair of pants (that is, a sphere from which three disjoint open discs have been removed) is  free on two generators. Fix a pair of generators, $a$ and $b$,  of this group,  each of them represented by a simple closed curve parallel to one of the  boundary components. Then every free homotopy class of closed curves  on the pair of pants can be represented uniquely as a  reduced cyclic  word in the symbols $a, b, \A, \B$. (A {\em cyclic word } $w$ is an equivalence class of words up to a cyclic permutation of their  letters; {\em Reduced}
means that the cyclic word contains no  juxtapositions of $a$ with $\A$,  or $b$ with $\B$.)
 The {\em word length $\wl(w)$} (with respect to the generating set $(a,b)$) of 
a free homotopy class of curves is the total number of letters occurring in
the corresponding reduced cyclic word.

 Consider a hyperbolic metric (that is, a metric with constant negative one curvature) on the pair of pants such that each of the  boundary components is a geodesic.
It is well known (see, for instance, \cite[Theorem 3.1.7]{bu}) that such hyperbolic metric on the pair of pants is determined by the lengths,  $A$, $B$ and $C$, of the three geodesics at each of the boundary components. After fixing $A$, $B$, and $C$ every free homotopy class of curves $w$ can be a	associated with a \emph{geometric length $G\abc(w)$}, the length of the unique geodesic in that class. 

The goal of this note is to study how the geometric length $\gl(w)$ varies  over the population $\FF{L}$ of all reduced cyclic words $w$ of  word length $L$ for each   positive integer $L$. (The set $\FF{L}$  is endowed with the uniform measure)

We choose a few representative metrics (mainly, the ratio between boundary components is about the same or very different). Using the parametrization described in Section \ref{Section: param},  we computed the length of a sample of $100,000$ words of $100$ letters and studied the distribution.

\begin{figure}[htp]
\includegraphics[width=.45\textwidth]{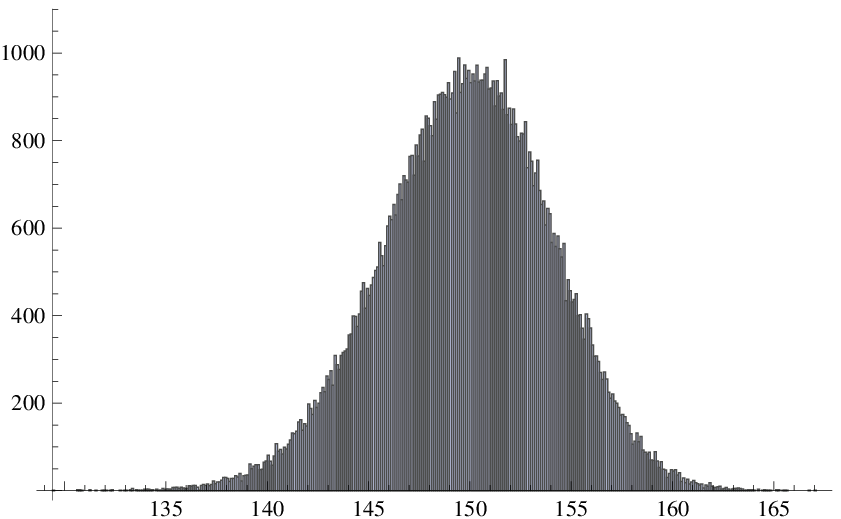}
\includegraphics[width=.45\textwidth]{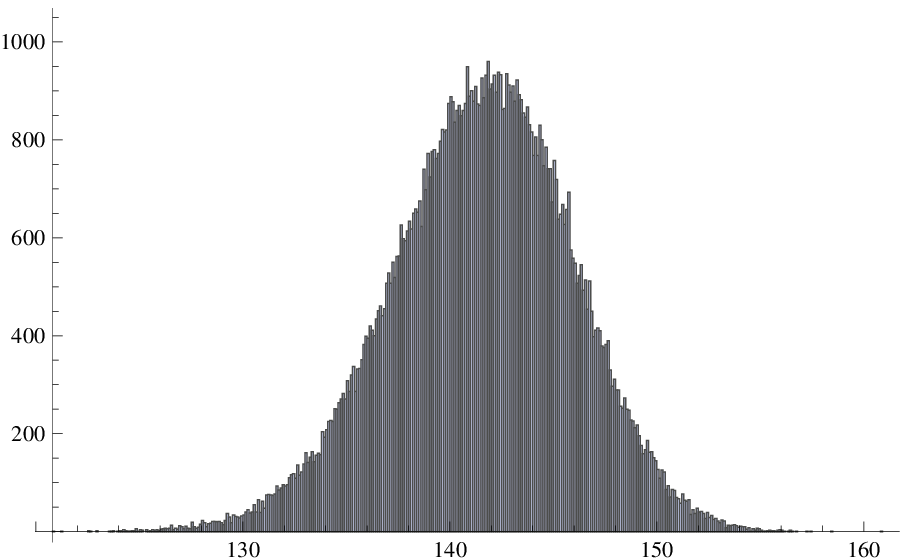}\\

\includegraphics[width=.45\textwidth]{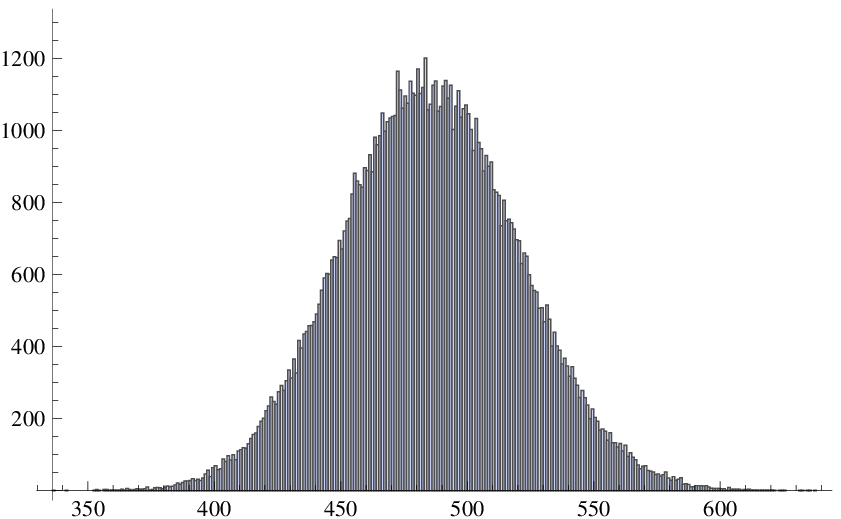}
\includegraphics[width=.45\textwidth]{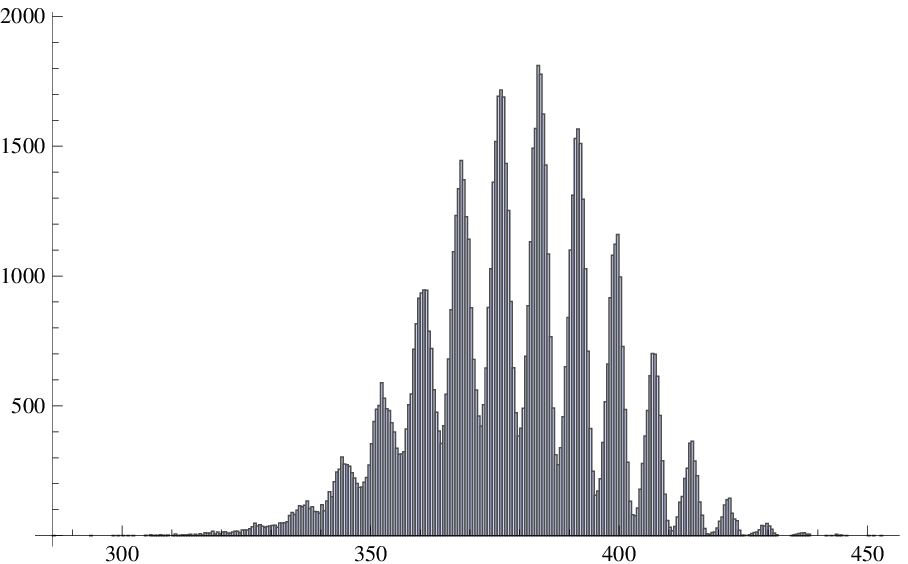}

\caption{Histograms of the geometric length of a sample of 100,000 words of word length 100.  The parameters are $(A,B,C)$;  $(1,1,1)$  top left, $(0.1, 1,1)$ top right, $(1,10,0.1)$  bottom right; $(0.1,1,10)$ bottom left}
\label{Fig: normal histogram}
\end{figure}

 The outcome of our experiments lead us to formulate the next conjecture. (see Figure \ref{Fig: normal histogram}). The data on which the histograms are based can be found  \cite{data}.

\begin{proclama}{Conjecture} For each hyperbolic metric  on the pair of pants with geodesic boundary and such that that length of the boundary components are $A$, $B$ and $C$, there exist constants $\kappa=\kappa(A,B,C)$ and $\sigma=\sigma(A,B,C)$ such that if a reduced cyclic word is chosen at random from among all classes in $\FF{L}$, then for large $L$ the distribution of the geometric length  approaches the Gaussian distribution with mean $\kappa \cdot L$ and standard deviation $\sigma \cdot L$.

More precisely, there exists constants $\kappa=\kappa(A,B,C)$ and $\sigma = \sigma(A,B,C)$, such that  for any $a < b$ the proportion of words $w$ in $\FF{L}$ such that 
$$
\frac{G\abc (w) - \kappa L}{\sqrt{L}} \mbox{ in } [a,b]
$$
converges, as $L$ goes to infinity to $\frac{1}{\sigma \sqrt{2\pi}}\int^b_a \exp\{-x^2/2\sigma^2\}dx$.
\end{proclama}

\begin{remark} 
The histogram of the sample $(0.1,1,10)$ in Figure \ref{Fig: normal histogram} does not suggest Gaussian distribution. The "peaks" in the graph are due to the following: The length of $C$ of one of the boundary components is very large compared to $A$ and $B$. Therefore, the length of a cyclic reduced word $w$ will be roughly, $C/2$ times the number of occurrences of $ab, ba, a\B, \B a, \A b, b\A, \A \B$ and $\B \A$ in $w$. 

Since it is not computationally feasible to study very long words with the metric $(0.1,1,10)$, we study instead the metric $(1,1,5)$. When $L$ is  small compared with $C$ ($14$ and $20$), we still see "peaks" in the histogram, approximately at multiples of $C/2=2.5$. The peaks start to disappear when $L$ is $50$ and do not appear at all when $L$ is $100.$
 This suggests that the "peaks" will gradually disappear when the word length goes to infinity. 

On the other hand, difference between the histograms of the metrics $(1,10,0.1)$ and $(0.1, 1,10)$ comes from the fact that our choice of basis for the fundamental group of the pair of pants, although natural, is not symmetrical, in the sense that we arbitrarily chose two of the three boundary components.  
\end{remark}

\begin{figure}[htp]
\includegraphics[width=.45\textwidth]{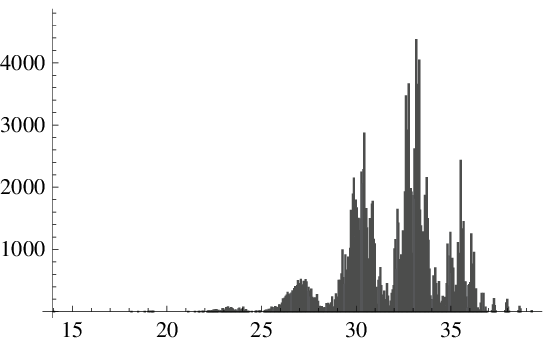}
\includegraphics[width=.45\textwidth]{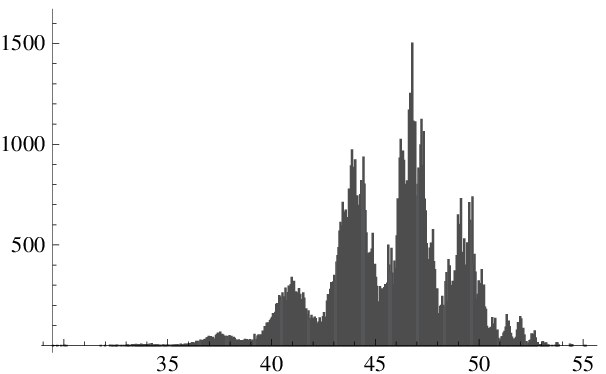}\\
\includegraphics[width=.45\textwidth]{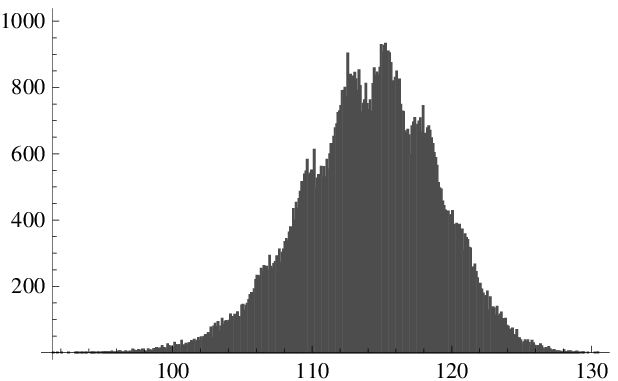}
\includegraphics[width=.45\textwidth]{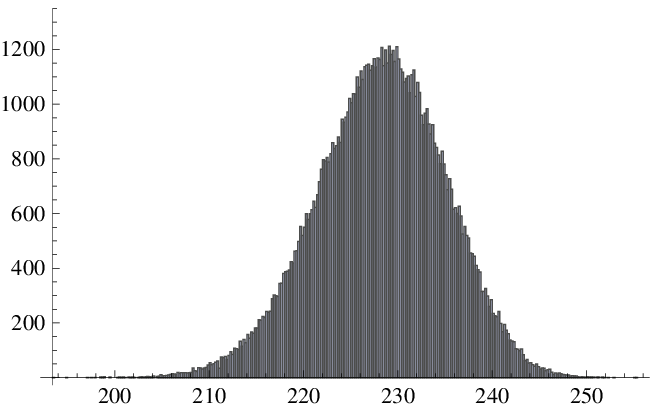}

\caption{ Top left: Histogram of all words of word length 14, with metric $(1,1,5)$. Top right, bottom left and bottom right respectively, are 
histograms of the geometric length of a sample of 100,000 words with parameters $(1,1,5)$ and word length $20$, $50$ and $100$ respectively.
}
\label{Fig: sequence}
\end{figure}

This type of experimental study was initiated  by Chas \cite{chas} who computed the distribution of the self-intersection number of free homotopy classes, sampling by word length. Chas and Lalley \cite{cl}   proved that  this distribution, properly normalized, approaches a Gaussian distribution when the word length goes to infinity.

One can also think of our experiments as a certain product of random matrices.

During the past decades there have been many results about the limiting behavior of products of random matrices. Bellman \cite{bellman} is the first to examine this type of question.  Fustebenger \cite{F} andFurstenberg and Kesten \cite{fk} stated a central limit theorem under certain hypothesis. Recently, Pollicott and Sharp \cite{ps} studied central limit theorems and their generalizations for matrix groups acting co-compactly or convex co-compactly on the hyperbolic plane. It is likely that their results, about infinite products, can be adapted to prove our conjecture, which concerns finite products.

\section{The parametrization of the hyperbolic metric on the pair of pants}\label{Section: param}

Let $a$ and $b$ be generators as above of the fundamental group $G$ of a pair of pants that has been endowed with a hyperbolic metric. We can regard $G$ as a subgroup of $PSL(2,{\mathbb R})$, acting on the upper half plane. There is then a unique correspondence between closed geodesics on the pair of pants and non-trivial elements of $G$, where the length of the geodesic, $l(w)$ is related to the absolute value of the trace of the corresponding element $g \in G$, by

$$\cosh \left(\frac{l}{2}\right)=\frac{|\rm{tr}(M)|}{2}.$$

The group $G$ is defined by the hyperbolic metric only up to conjugation in $PSL(2,{\mathbb R})$.  We normalize $G$ so that first, the line orthogonal to the axes of $a$ and $b$ is the imaginary axis, where $0$ is closer to the axis of $a$ than to the axis of $b$. If necessary we replace $a$ and/or $b$ by its inverse so that $a(0) > 0$ and $b(0) < 0)$. We then normalize further so that $a$ can be realized as reflection in the imaginary axis followed by reflection in the hyperbolic line whose endpoints are at $\alpha$ and $1$, $\alpha < 1$. With this normalization, $b$ can be realized as reflection in the hyperbolic line whose endpoints are at $\beta$ and $\gamma$, $\beta < \gamma$, followed by reflection in the imaginary axis. 

Since the product $a\cdot b$ is necessarily hyperbolic, we note that we have the inequalities, 
\begin{equation}\label{eqn:ineq}
0 < \alpha < 1 < \beta < \gamma.
\end{equation}

Conversely, if $\alpha$, $\beta$ $\gamma$ satisfy the above inequality, and $G$ is the group generated by $a$ and $b$, where $a$ is defined as reflection in the imaginary axis followed by reflection in the line with endpoints $\alpha$ and $1$, and $b$ is defined as reflection in the line whose endpoints are at $\beta$ and $\gamma$, followed by reflection in the imaginary axis, then $G$ is necessarily is a Fuchsian group representing a pair of pants, and $a$ and $b$ represent simple boundary geodesics. Easy computations show that we now have

$$ \hat a = \frac{1}{1-\alpha}\left(
\begin{array}{cc}
1+\alpha & 2\alpha \\
2& 1+\alpha \\
\end{array} \right),\ \
\hat b = \frac{1}{\gamma-\beta}\left( \begin{array}{cc}
\beta+\gamma & -2\beta\gamma  \\
-2 & \beta+\gamma \\
\end{array} \right).$$

As above, let $A$, $B$ and $C$ be the lengths of the geodesics corresponding to $a$, $b$, and $c = ab$; these are the lengths of the three boundary geodesics on the pair of pants.

Then
$$x = \cosh \frac{A}{2} = \frac{|\rm{tr}(a)|}{2} = \frac{1+\alpha}{1-\alpha},$$
$$y = \cosh\frac{B}{2} = \frac{|\rm{tr}(b)|}{2} = \frac{\gamma + \beta}{\gamma - \beta}$$
and
$$z = \cosh\frac{C}{2} = \frac{\rm{tr}(ab)}{2} = -\frac{(1 + \alpha)(\beta + \gamma) - 2(\alpha + \beta\gamma)}{(1 - \alpha)(\gamma - \beta)}.$$

Using inequality \ref{eqn:ineq}, the above can be uniquely solved for $\alpha, \beta$ and $\gamma$. We obtain

$$\alpha=\frac{x-1}{x+1},$$
$$\beta=\frac{(xy+z)+\sqrt{x^2+y^2+z^2+2xyz-1}}{(x+1)(y+1)},$$
$$\gamma=\frac{(xy+z)+\sqrt{x^2+y^2+z^2+2xyz-1}}{(x+1)(y-1)}.$$

\end{document}